\input amstex
\documentstyle{amsppt}
\pagewidth{6.5truein}
\pageheight{9truein}
\TagsOnRight

\def\({\left(}
\def\){\right)}
\def\[{\left[}
\def\]{\right]}
\def\l|{\left|\right.}
\def\r|{\left.\right|}
\def\ov{\overline}
\def\noi{\noindent}
\def\ve{\varepsilon}
\def\wideh{\widehat}
\def\widet{\widetilde}

\def\vars{\varsigma}

\def\disp{\displaystyle}
\def\var{\varphi}
\def\cirdot{\circ\cdots\circ}

\topmatter
\title
The trace problem for vector fields\\
satisfying H\"ormander's condition.
\endtitle
\author
S. Berhanu and I. Pesenson
\endauthor

\address
Department of Mathematics, Temple University,
Philadelphia, PA 19122, U.S.A.
\endaddress
\email
berhanu\@euclid.math.temple.edu
\endemail

\address
Department of Mathematics, Temple University,
Philadelphia, PA 19122, U.S.A.
\endaddress
\email
pesenson\@euclid.math.temple.edu
\endemail

\define\eigen#1{\Cal E\(L_#1\)}
\define\auto#1#2{\Cal H_#1\(#2\)}
\define\part{\partial}

\abstract
Trace theorems are proved for non-isotropic Sobolev and
$L^p$-Lipschitz spaces
defined by vector fields satisfying H\"ormander's bracket
condition of order 2. It is shown
that the loss of regularity by traces is the same as in the
classical case.
\endabstract

\endtopmatter

\document
\baselineskip=16pt

\medskip
\heading
0. Introduction.
\endheading

It is a classical fact that for Sobolev spaces $W^r_p(\Bbb
R^{n+1})$, the space
of traces is the\break $L^p$-Lipschitz (Besov) space
$\land^{r-1/p}_p(\Bbb R^k)$.
This result when $p=2$ was obtained in \cite{1} and \cite{9},
and for $r=1$ and
$1<p<\infty$
in \cite{3}. The complete solution for the integer and
fractional Sobolev spaces was
obtained by E. Stein \cite{10} and for the Besov spaces by O.
Besov \cite{2}.

In this paper we consider the analogous problem for
non-isotropic Sobolev and
$L^p$-Lipschitz spaces. It turns out that as in the classical
case, the space of
all traces can be described in terms of some kind of Besov
norm constructed by
means of ``tangential components" of the given vector fields
and their one-parameter
groups of diffeomorphisms. We have the same phenomenon as in
the classical situation:
traces are less regular than the original functions and the
loss of regularity
is precisely $\disp{1\over p}$. We prove both restriction and
extension theorems that
are compatible. The extension result is established by
analyzing an explicit
extension operator which is a non-isotropic version of the
classical Hardy operator.  In the model case on the Heisenberg
group, this problem was solved in \cite{5} and \cite {6}.

This article is organized as follows. In section 1 we first
state our results and
prove the independence of our function spaces on the bases
used. We then present
the proof of our restriction result, Theorem 1.1. In section
2 we present our
extension operator and prove the extension theorem, Theorem
1.2.
\bigskip
\bigskip

\heading
1. Statement of results and independence of bases.
\endheading

For a point in ${\Bbb R}^{n+1}$, we will use coordinates $(x,
t)$ where $x\in {\Bbb
R}^n$ and $t\in {\Bbb R}$ and view ${\Bbb R}^{n+1}$ as ${\Bbb
R}^n_x\times {\Bbb
R}_t$. We will also identify the subset ${\Bbb R}_x^n\times
\{0\}$ with ${\Bbb
R}_x^n$. For any $y \in {\Bbb R}^{n+1}$, the vector space
$T_y{\Bbb R}_x^n$ will
denote
$$\left \{\sum_{j=1}^n a_j\frac {\partial}{\partial x_j} :
a_j\in {\Bbb R}\text{ for }
j = 1,\dots, n\right \}.$$
Let ${\Cal V}$ be a $C^{\infty}$ vector subbundle of the
tangent space $T{\Bbb
R}^{n+1}$ near $0$. Let the fiber dimension of ${\Cal V}$ be
$k+1$. For any point $y$
where ${\Cal V}$ is defined, ${{\Cal V}}_y$ will denote the
fiber of ${\Cal V}$ at
$y$. We will assume that ${\Cal V}$ satisfies the following
two conditions:

\medskip
(i) ${\Cal V}_0\subsetneq T_0\Bbb R^n_x$,\qquad and

(ii) The sections of ${\Cal V}\cap T\Bbb R^n_x$ together with
their brackets $[X,Y]$
span $T\Bbb R^n_x$ near $0$ in $\Bbb R^n_x$.
\medskip
Assumption (i) means that there is a vector $v$ in the fiber
${{\Cal V}}_0$ with a
nonzero $\frac {\partial}{\partial t}$ component. It follows
that ${\Cal V}\cap T{\Bbb
R}_x^n$ forms a bundle of fiber dimension $k$ near $0$.
Indeed, since ${{\Cal V}}_0 +
T_{0}{\Bbb R}_x^n = T_{0} {\Bbb R}^{n+1}$, by continuity,
${{\Cal V}}_y + T_{y}{\Bbb
R}_x^n = T_{y} {\Bbb R}^{n+1}$ for $y$ near $0$. Hence
${{\Cal V}}_y \cap T_{y}{\Bbb
R}_x^n$ is of dimension $k$ for $y$ near $0$, telling us that
${\Cal V}\cap T{\Bbb
R}_x^n$ is a bundle. Condition (ii) therefore says that the
restriction of this bundle
to ${\Bbb R}^n_x\times \{0\}$ satisfies H\"ormander's bracket
condition of order 2.

\medskip
Here is a simple example in ${\Bbb R}^4 = {\Bbb R}^3_x\times
{\Bbb R}_t$, where
$x=(x_1,x_2,x_3)$. Let ${{\Cal V}}^{\prime}$ be the
$C^{\infty }$ bundle generated by
$\frac {\partial}{\partial x_1},\frac {\partial}{\partial
t},\text{ and } \frac
{\partial}
{\partial x_2}+x_1\frac {\partial}{\partial x_3}$. Then
$\frac {\partial}{\partial
t}\notin T_0{\Bbb R}^3_x$ and so (i) is met. Since ${{\Cal
V}}^{\prime}\cap T{\Bbb
R}^3_x$ is generated by $\frac {\partial}{\partial x_1}\text{
and } \frac {\partial}
{\partial x_2}+x_1\frac {\partial}{\partial x_3},\text{ and }
[\frac
{\partial}{\partial x_1}, \frac {\partial}{\partial
x_2}+x_1\frac {\partial}{\partial
x_3}]=\frac {\partial}{\partial x_3},$ we see that (ii) is
also met.

\medskip
Let $\beta=\{Z_1,\cdots,Z_k\}$ be a basis of ${\Cal V}\cap
T\Bbb R^n_x$ over an open
neighborhood $V$ of $0$ in $\Bbb R^n_x$. Let $V_1$ be a
neighborhood of $0$
such that $V_1\subset\subset V$ and suppose $\delta >0$
satisfies

$$
e^{\tau Z_j}(V_1)\subseteq V ~~\hbox{ for } |\tau|\le\delta
~~\hbox{ and for
all } j.
$$
(Here $e^{\tau Z_j}x$ denotes the integral curve of $Z_j$
starting at $x$ when
$\tau=0$). Let $1<p<\infty$. For $\psi\in C^\infty_0(V_1)$,
define

$$
\omega(t,\psi,Z_j,V_1,V)=\sup_{|\tau|\le t} \|e^{\tau
Z_j}\psi-\psi\|_{L^p}
$$
and

$$
\|\psi\|_{W_{1-{1\over
p},p}(\beta,V,V_1,\delta)}=\|\psi\|_{L^p}+\sum^k_{i=1}
\left\{\int^\delta_0\[t^{-\theta}\omega_i(t,\psi,V)\]^p{dt\ov
er t}\right\}^{1
\over p},
$$
where $\theta=1-{1\over p}$ and
$\omega_i(t,\psi,V)=\omega(t,\psi,Z_i,V_1,V)$.

Note that if $0<\delta'<\delta$, then

$$
\|\psi\|_{W_{1-{1\over p},p}(\beta,V,V_1,\delta)}
$$
is equivalent to

$$
\|\psi\|_{W_{1-{1\over p},p}(\beta,V,V_1,\delta')}
$$
and hence, in the sequel, we'll simply write

$$
\|\psi\|_{W_{1-{1\over p},p}(\beta,V,V_1)}
$$
with the implicit understanding that we are using some
$\delta >0$ satisfying

$$
e^{\tau Z_j}(V_1)\subseteq V \hbox{ for } |\tau|\le\delta
\hbox{ and } j=1,\cdots,k.
$$

We will next show that if we change the basis $\beta$ to
$\beta'$, then after
contracting the neighborhoods $V$ and $V_1$, the norms become
equivalent.
More precisely, we have:

\proclaim{Lemma 1.1} Let $\beta=\{Z_1,\cdots,Z_k\}$ and
$\beta'=\{Y_1,\cdots,Y_k\}$
be bases of ${\Cal V}\cap T\Bbb R^n_x$ over a neighborhood
$V$ of $0$ in $\Bbb R^n_x$.
Then there exist neighborhoods $V_2\subset\subset
V_1\subset\subset V$ and
$C>0$ such that for all $\psi\in C^\infty_0(V_2)$,

$$
\|\psi\|_{W_{1-{1\over p},p}(\beta',V_2,V_1)}\le
C\|\psi\|_{W_{1-{1\over p},p}
(\beta,V_2,V_1)}.
$$
\endproclaim

\demo{Proof} Since the $Z_j$ together with their brackets
span $T\Bbb R^n_x$
near $0$, after contracting $V$ if necessary, we get a basis

$$
\widet\beta=\{Z_1,\cdots,Z_k,Z_{k+1},\cdots,Z_n\}
$$
of $T\Bbb R^n_x$ over $V$ where for each $i\ge k+1$,
$Z_i=\[Z_1^i,Z_2^i\]$ for
some $Z^i_1,Z^i_2$ in $\beta$.

Let $Y\in \beta'$. Write $Y=\sum^k_{j=1} a_j(x)Z_j$ for some
$a_j\in C^\infty
(V)$. To prove the Lemma, we need to dominate

$$
\(\int\left|\psi\(e^{\tau Y}x\)-\psi(x)\right|^p dx\)^{1\over
p}=\(\int\left|
\psi\(e^{\tau(a_1(x)Z_1+\cdots+a_k(x)Z_k)}x\)-\psi(x)\right|^
p dx\)^{1\over p}
$$
by terms of the form

$$
\(\int|\psi(e^{sZ_j}x)-\psi(x)|^p dx\)^{1\over p}
$$

For each $k+1\le i\le n$, we will define mappings $F_i(s)(x)$
which are approximations
of $e^{sZ_i}x$.

Recall that for such $i$, $Z_i=[Z^i_1,Z^i_2]$ where $Z^i_1$
and $Z^i_2$ are in
$\beta$.

Define
\medskip

$$
F_i(s)(y)=\cases e^{-\sqrt{s}Z^i_2} e^{-\sqrt{s}Z^i_1}
e^{\sqrt{s}Z^i_2} e^{
\sqrt{s}Z^i_1}(y), & s\ge 0\\
e^{-\sqrt{|s|}Z^i_1} e^{-\sqrt{|s|}Z^i_2} e^{\sqrt{|s|}Z^i_1}
e^{\sqrt{|s|}Z^i_2}(y),
& s< 0. \endcases
$$
\medskip

By the Campbell-Hausdorff formula,

$$
F_i(s)(x) = e^{sZ_i}g(x,s)
$$
where

$$
g(x,s)=x+0\(|s|^{3/2}\)\tag $\ast$
$$
and hence each $F_i$ is $C^1$. For $s=(s_1,\cdots,s_n)$ and $x$ small, define
 $$F(s,x)=e^{s_1Z_1}\cdots e^{s_kZ_k}F_{k+1}(s_{k+1})\cdots F_n(s_n)x$$

The estimate in $(\ast)$ tells us that for each $x$, $F(s,x)$
is a $C^1$
diffeomorphism from a neighborhood of $0$ in $s$ space to a
neighborhood of $x$.
In fact, there is $\ve >0$ and neighborhoods $V_1$ and $V_2$
of $0$,
$V_2\subset\subset V_1$, such that for each $x$ in $V_2$,
$s\longmapsto F(s,x)$
is a diffeomorphism from $B_\ve (0)$ into $V_1$.

It follows that for $\tau$ near $0$, the implicit function
theorem gives us
functions

$$
s(\tau,x)=(s_1(\tau,x),\cdots,s_n(\tau,x))
$$
such that $s(\tau,x)=o(|\tau|)$ and

$$
F(s(\tau,x),x)=e^{\tau(a_1(x)Z_1+\cdots+a_k(x)Z_k)}x=e^{\tau
Y} x.
$$

We therefore need to dominate terms of the form

$$
\(\int\left|\psi\(e^{s_1(\tau,x)Z_1}\cdots
e^{s_k(\tau,x)Z_k}\cdots
e^{b_1(\tau,x)X_1}\cdots
e^{b_\ell(\tau,x)X_\ell}x\)-\psi(x)\right|^p dx\)^{1\over p},
$$
where the $X_i\in \beta$, $b_j(\tau,x)=o(|\tau|)$,

$$
s_i(\tau,x)=o(|\tau|).
$$

After using the triangle inequality and change of variables,
we are led to terms
of the form

$$
\(\int\left|\psi\(e^{b(\tau,x)Z}x\)-\psi(x)\right|^p
dx\)^{1\over p},
$$
where $Z\in \beta$ and $b(\tau,x)=o(|\tau|)$. Finally, an
application of the
technique used to prove Lemma 4.1 in \cite{4} enables us to
dominate these latter
terms by integrals of the form

$$
\sup_{|s|\le
C|\tau|}\(\int\left|\psi\(e^{sZ}x\)-\psi(x)\right|^p
dx\)^{1\over p},
$$
where $C$ is independent of $\psi$. The lemma follows from
these observations.

If $\beta=\{Y_1,\cdots,Y_{k+1}\}$ is a basis of ${\Cal V}$
over a neighborhood  $U$ of
$0$ in $\Bbb R^{n+1}$, and $\var\in C^\infty_0 (U)$, we
define

$$
\|\var\|_{W_{1,p}(U,\beta)}=\|\var\|_{L^p}+\sum^{k+1}_{j=1}
\|Y_j\var
\|_{L^p}.
$$

It is clear that if $\widetilde\beta$ is also a basis of
${\Cal V}$ over $U$, we get
an equivalent norm. Hence in the sequel, we'll often omit
mention of the basis.

We are now ready to state the main results of this article:
\enddemo

\proclaim{Theorem 1.1} Let $1<p<\infty$. Let $\beta$ and
$\beta'$ be any bases
near $0$ of ${\Cal V}$ and ${\Cal V}\cap
T\Bbb R^n_x$ respectively.
Then there exist neighborhoods $U$ of $0$ in $\Bbb R^{n+1}$
and $V$ of $0$ in
$\Bbb R^n_x$ and $C>0$ such that if $\var\in C^\infty_0(U)$
and $R\var(x)
=\var(x,0)$, then

$$
\|R\var\|_{W_{1-{1\over p},p}(V,\beta')}\le
C\|\var\|_{W_{1,p}(U,\beta)}.
$$
\endproclaim

Conversely, we'll prove the following extension theorem.

\proclaim{Theorem 1.2} Let $1<p<\infty$. Let $\beta$ and
$\beta'$ be any bases
near $0$ of ${\Cal V}$ and ${\Cal V}\cap
T\Bbb R^n_x$ respectively.
Then there exist neighborhoods $U$ of $0$ in $\Bbb R^{n+1}$
and $V$ of $0$ in
$\Bbb R^n_x$ and a linear extension mapping from
$W_{1-{1\over p},p}(V,\beta')$
to $W_{1,p}(U,\beta)$ that is continuous.
\endproclaim
\bigskip

\noindent{\bf Remark 1.}~~ Theorem 1.2 shows that the loss
${1\over p}$ of smoothness
in Theorem 1.1 is sharp.
\bigskip
\noindent{\bf Remark 2.}~~ As indicated in the introduction,
these theorems show
that traces lose exactly the same smoothness as in the
classical case.
\bigskip

\demo{Proof of Theorem 1.1} We begin by observing that we can
choose sections
$X_1,\cdots,X_k$ of ${\Cal V}$ of the form

$$
X_i=\sum^n_{j=1} a_{ij}(x,t){\part\over\part x_j}, \qquad
1\le i\le k
$$
such that $\beta'=\{Y_1,\cdots,Y_k\}$ where
$Y_i=X_i\left|_{t=0}\right.$

To see this, let $\beta'=\{Y_1,\cdots,Y_k\}$ and choose a
basis
$Z=\{Z_1,\cdots,Z_{k+1}\}$
of ${\Cal V}$ of the form
$$
Z_i=\sum^n_{j=1} b_{ij}(x,t){\part\over\part x_j}\quad\hbox{
for } 1\le i
\le k
$$
and

$$
Z_{k+1}={\part\over\part t}+\sum^n_{j=1}
C_j(x,t){\part\over\part
x_j}.
$$

Such a basis $Z$ is possible  since ${\Cal V}$ is not
contained in $T\Bbb R^n_x$.
Let $f_{ij}(x)$ be $C^\infty$ functions such that

$$
Y_i=\sum^n_{j=1}f_{ij}(x)Z_j|_{t=0}\quad\hbox{ for } 1\le
i\le k.
$$

Set $X_i=\sum^n_{j=1}f_{ij}(x)Z_j$ for $1\le i\le k$. Then
$\{X_1,\cdots,X_k\}$ is
as desired. Moreover, if $X_{k+1}=Z_{k+1}$, then
$\{X_1,\cdots,X_{k+1}\}$ is a
basis of ${\Cal V}$ near $0$.

Next we observe that we may assume $X_{k+1}$ to be
$\disp{\part\over\part t}$.

Indeed, suppose $G(x,t)$ is a diffeomorphism from $(x,t)$
space to $(y,s)$ space
such that $G(x,0)=(x,0)$ and
$G_\ast\(X_{k+1}\)=\disp{\part\over\part s}$.
Since $G(x,0)=(x,0)$, we observe that it suffices to prove
the theorem in
$(y,s)$ space for the bundles $G_\ast({\Cal V})$ and
$G_\ast({\Cal V})\cap T\Bbb
R^n_y$.
Thus we will assume that $\{X_1,\cdots,X_{k+1}\}$ is a basis
of ${\Cal V}$ near $0$,
$X_i|_{t=0}=Y_i$ for $1\le i\le k$ and
$X_{k+1}=\disp{\part\over\part t}$.

Let $U=V\times (-\ve,\ve)$ be a neighborhood of $0$ in $\Bbb
R^
{n+1}$ over which the $X_j$ span ${\Cal V}$. Fix
$X\in\{X_1,\cdots,X_k\}$, and let
$L=\disp{\part\over\part t}-X$. If $\var(x,t)\in
C^\infty_0(U)$, we will
express $\var(x,t)$ in terms of $L\var(x,t)=f(x,t)$ and
$\var_0(x)=
\var(x,0)$ as follows.

Let $p_j(x,t)$ $(1\le j\le n)$ be the unique solution of
\medskip

$$
\cases Lp_j(x,t) &= 0 \\
p_j(x,0) & =x_j \endcases
$$
\medskip
\noindent in a neighborhood of $0$ which we still call $U$.
Define $G(x,t)=
\(p(x,t),t\)$ where $p=(p_1,\cdots,p_n)$. Consider the change
of variables

$$
(x,t)\longmapsto (y,s)=(p(x,t),t).
$$

If $g=g(y,s)$, we have:

$$
L(g(G(x,t)))={\part g\over \part s}(G(x,t)).
$$
Hence if $H(y,s)=(h(y,s),s)$ is the inverse of $G$ and  $F$
solves $\disp{\part F
\over\part s}(y,s)=f(H(y,s))$, $F(y,0)=\var_0(y)$, then

$$
\var(x,t)=F(G(x,t)).
$$
Hence

$$
\var(x,t)=\var_0(p(x,t))+\int^t_0 f(H(p(x,t),\tau)) d\tau
\tag 1.1
$$

If $\psi=\psi(x)$, let $G(t)\psi(x)$ denote the function
$\psi(p(x,t))$. Using
this notation we can write

$$
\var(x,t)=G(t)\var_0(x) +\int^t_0 f(H(p(x,t),\tau)) d\tau
\tag 1.1'
$$

The proof of Theorem 1.1 will use the following: ($X$ will
continue to denote
an element of $\{X_1,\cdots,X_k\}$.)

\enddemo

\proclaim{Lemma 1.3} There exist neighborhoods $V,U$ of $0$
in $\Bbb R^n$ and
$\Bbb R^{n+1}$ respectively, $\delta >0$ and $C>0$ such that
for any $\var\in
C^\infty_0(U)$,

$$
\left\{\int^\delta_0\[t^{{1\over p}-1}\sup_{|\tau|\le
t}\|G(\tau)\var_0-\var_0\|
_{L^p(V)}\]^p {dt\over t}\right\}^{1\over p}\le
C\(\left\|{\part\var\over
\part t}\right\|_{L^p(U)}+\left\|X\var\right\|_{L^p(U)}\).
$$

\endproclaim

\demo{Proof of Lemma 1.3} We take $U=V\times(-\ve,\ve)$ so
that (1.1') is valid.
From (1.1') we have:

$$
G(\tau)\var_0(x)-\var_0(x)=\int^\tau_0{\part\var\over\part
s}(x,s) ds-
\int^\tau_0 f(H(p(x,\tau),s)) ds.
$$

Since the $x$ support of $f(x,t)=L\var(x,t)$ is in $V$,
Minkowski's inequality
yields

$$
\|G(\tau)\var_0-\var_0\|_{L^p(V)}\le
C\left\{\int^\tau_0\left\|{\part\var
\over\part s}(\cdot,s)\right\|_{L^p(V)}
ds+\int^\tau_0\|f(\cdot,
s)\|_{L^p(V)} ds\right\}.
$$

Thus for any $t\in [0,\ve)$, we have:

$$
t^{-1}\sup_{|\tau|\le t}\|G(\tau)\var_0-\var_0\|_{L^p(V)}\le
C\left\{t^{-1}\!\!
\int^t_0\left\|{\part\var
\over\part s}(\cdot,s)\right\|_{L^p(V)}
ds+t^{-1}\!\!\int^t_0\|L\var
(\cdot,s)\|_{L^p(V)} ds\right\}.
$$

To the latter we apply the Hardy-Littlewood inequality to get
the Lemma for
any $\delta\le\ve$.

\noi(Recall that the Hardy-Littlewood inequality says that

$$
\left\{\int^\infty_0\left| t^{-1}\int^t_0 h(s)ds\right|^q
dt\right\}^{1\over q}
\le C\(\int^\infty_0\left|h(\tau)\right|^q d\tau\)^{1\over q}
$$
for $1<q<\infty$).
\enddemo

End of the proof Theorem 1.1. As indicated already, we may
let

$$
X_i=\sum^n_{j=1} a_{ij}(x,t){\part\over\part x_j}\quad\hbox{
for }\quad
1\le i\le k,
$$

$$
X_i\big|_{t=0}=Y_i\qquad\qquad\hbox{and}\qquad\qquad
X_{k+1}={\part\over\part
t}.
$$

For each $i=1,\cdots,k$, let $B_i(x,t)=e^{tY_i}x$ where we
view $Y_i$ as a vector
field in $\Bbb R^{n+1}$. Let $p^i=\(p^i_1,\cdots,p^i_n\)$ for
$i\le i\le k$ be
the unique solution of
\medskip

$$
\cases \disp{\part p^i_j\over \part t}(x,t)-X_ip^i_j &=0\\
p^i_j(x,0) &=x_j. \endcases
$$
\medskip

Since $p^i(x,0)=B_i(x,0)$ and $X_i\big|_{t=0}=Y_i$, we have:

$$
p^i(x,t)=B_i(x,t) +0(t^2).
$$

Therefore, if $R_i(t)x$ denotes $e^{-tY_i}p^i(x,t)$, then

$$
R_i(t)x=x+o(t^2)\qquad\hbox{ and }\qquad
R_i(t)^{-1}x=x+o(t^2).
$$

Let $R^i(t)=R_i(t)^{-1}$. We have:

$$
e^{tY_i}y=p^i\(R^i(t)y,t\)=G_i(t)(R^i(t)y)
$$
where we have used the notation $G_i(t)x=p^i(x,t)$.

Let $V'$ be a neighborhood of $0$ such that $V'\subseteq V$
and

$$
R^i(\tau)(V')\subseteq V\quad\forall~i=1,\cdots,k\hbox{ and
for } 0\le\tau\le
\delta_1,~\delta_1<\ve.
$$

Let

$$
\omega_i(t,\var_0,V')=\sup_{|\tau|\le t}\left\|e^{\tau
Y_i}\var_0-\var_0\right\|
_{L^p(V')}.
$$

From

$$
\var_0\(e^{\tau
Y_i}x\)-\var_0(x)=\var_0\(G_i(\tau)\(R^i(\tau)x\)\)-\var_0\(R
^i
(\tau)x\)+\var_0\(R^i(\tau)x\)-\var_0(x),
$$
we have, for $0\le t\le\delta_1$,

$$
\omega_i(t,\var_0,V')\le C\left\{\sup_{|\tau|\le
t}\left\|R^i(\tau)\var_0-\var_0
\right\|_{L^p(V)}+\sup_{|\tau|\le
t}\left\|G_i(\tau)\var_0-\var_0\right\|_{L^p
(V)}\right\}\tag 1.2
$$

Recall that $R^i(\tau)x=x+o(\tau^2)$ and so by Lemma 3.4 in
H\"ormander (\cite{4}),
we get:

$$
\sup_{|\tau|\le
t}\left\|R^i(\tau)\var_0-\var_0\right\|_{L^p(V)}\le
C~\omega(t^2,\var_0,V) \tag 1.3
$$
where

$$
\omega(t^2,\var_0,V)=\sup_{|s|\le t^2}\|\var(\cdot
+s)-\var(\cdot)\|_
{L^p(V)}
$$
is the usual $L^p$ modulus of continuity.

From the inequalities (1.2), (1.3) and Lemma 1.3 we get

$$
\align
& \left\{\int^{\delta_1}_0\[t^{{1\over
p}-1}\omega_i(t,\var_0,V')\]^p{dt\over
t}\right\}^{1\over p}\\
\le\quad & C\left\{\(\int^{\delta_1}_0 t^{{1\over
p}-1}\sup_{|\tau|\le t}
\left\|G_i(\tau)\var_0-\var_0\right\|^p_{L^p(V)}{dt\over
t}\)^{1\over p}\right.\\
&\quad\left. +\(\int^{\delta_1}_0\[t^{{1\over
p}-1}\omega(t^2,\var_0,V)\]^p
{dt\over t}\)^{1\over p}\right\}\\
\le\quad &
C\left\{\|\var\|_{W_{1,p}(U)}+\(\int^{\delta_1}_0\[t^{-\sigma
}\omega
(t,\var_0,V)\]^p{dt\over t}\)^{1\over p}\right\} \tag 1.4
\endalign
$$
where $\sigma={1\over 2}\(1-{1\over p}\)$.

The term
$\(\int^{\delta_1}_0\[t^{-\sigma}\omega(t,\var_0,V)\]^p{dt\ov
er t}\)^
{1\over p}$ is the main part of the norm in the Besov space
$B^\sigma_p(V)$.

We claim that for $V$, $U$ small enough, $\exists\, C>0$ such
that

$$
\|\var_0\|_{B^\sigma_p(V)}\le C\|\var\|_{W_{1,p}(U)}.\tag 1.5
$$

Indeed, first note that since $\left\{X_i\big|_{t=0} :1\le
i\le k\right\}$
satisfy the H\"ormander condition of order $2$, if the
neighborhood $U$ is
small enough, the fields $\{X_i : 1\le i\le k+1\}$ will
satisfy the same
condition in $U$. Hence by a result in \cite{7} we have

$$
\|\var\|_{L^{1\over 2}_p\(\Bbb R^{n+1}\)}\le
C\|\var\|_{W_{1,p}(U)}\quad\hbox{
for }\var\in C^\infty_0(U), \tag 1.6
$$
where $L^{1\over 2}_p \(\Bbb R^{n+1}\)$ is the space of
Bessel potentials in
$\Bbb R^{n+1}$ (see \cite{10} for definition).

Thus

$$
\left\|{\part\var\over\part t}\right\|_{L^p}
+\left\|\var\right\|_{L^{1\over
2}_p}\le C\|\var\|_{W_{1,p}(U)}\tag 1.7
$$
for $\var\in C^\infty_0(U)$ since
$X_{k+1}=\disp{\part\over\part t}$.

Next we recall the trace theorem (see \cite{11})

$$
\|\var_0\|_{B^\sigma_p(V)}\le C\(\left\|{\part\var\over\part
t}\right\|_
{L^p}+\|\var\|_{L^{1\over 2}_p}\). \tag 1.8
$$

From (1.7) and (1.8), we get

$$
\|\var_0\|_{B^\sigma_p(V)}\le C\|\var\|_{W_{1,p}(U)}.
$$

The latter together with inequality (1.4) prove the theorem.
\bigskip
{\bf Remark 3.} Since our vector fields satisfy H\"ormander's
bracket condition of
order 2, we were able to use inequalities $(1.6)$ and
$(1.8)$. Although the version of
(1.6) for commutators of all orders is known (see [5]), we
have not been able to
exploit it to get a reasonable generalization of Theorem 1.1.
\medskip
At this point, before we proceed to the extension theorem,
we would like to stress the implications of Lemma 1.3. This
lemma tells us that
even when the vector field $X$ is singular on $\Bbb R^n_x$,
the restriction
$\var_0(x)=\var(x,0)$ may gain some smoothness along some
direction in $x$ space.

As an example, let $X=t^m\disp{\part\over\part x_1}$ where
$m$ is a positive
integer. In the notation used in the lemma, we get
information on the $L^p$
modulus of

$$
\var_0\(x_1+ {\tau^{m+1}\over
m+1},x_2,\cdots,x_n\)-\var_0\(x_1,x_2,\cdots,x_n\).
$$
\bigskip
\bigskip

\heading
2. The extension operator.
\endheading
\bigskip

We now fix a special basis
$\left\{X_1,\cdots,X_k,\disp{\part\over\part t}\right
\}$ of ${\Cal V}$ of the form

$$
X_i={\part\over\part x_i}+\sum^n_{j=k+1}
a_{ij}(x,t){\part\over\part
x_j},\qquad 1\le i\le k
$$
which is achieveable after a permutation of the $x$
coordinates. The vector field
$\disp{\part\over \part t}$ comes after using a
diffeomorphism that preserves
the $x$ space as we saw before.

Let $Z_i=X_i\big|_{t=0}$ for $1\le i\le k$. By the hypotheses
on ${\Cal V}\cap T\Bbb
R^n_x$, the $Z_i$ together with their brackets span $T\Bbb
R^n_x$ near $0$.
We may therefore choose $Z_{k+1},\cdots,Z_n$ such that
$\{Z_1,\cdots,Z_n\}$ is a basis
of $T\Bbb R^n$ and for $i>k$ each $Z_i$ has the form $\[Z_\ell,Z_m\]$
for some $\ell,m
\le k$.

For $V$ a sufficiently small neighborhood of $0$ in $\Bbb
R^n_x$ and $t$ small,
define
\medskip

$$
H_i\var(x,t)=\cases \disp{1\over
t}\disp{\int^t_0}\var\(e^{\tau Z_i}x\)d\tau, & i\le
k\\
\disp{1\over t^2}\disp{\int^{t^2}_0}\var\(e^{\tau
Z_i}x\)d\tau, & i\ge k+1,\endcases
$$
\medskip
\noi where $t>0$ and $\var\in C^\infty_0(V)$.

Define $H\var(x,t)=(H_1\cirdot H_n\var)(x,t)$, where for
$\psi=\psi(x,t)$, $H_i
\psi(x,t)$ is defined by letting $H_i$ act on the function

$$
x\longmapsto \psi(x,t).
$$

For $\var\in C^\infty_0(V)$, define

$$
E\var(x,t)=\cases H\var(x,t), &t\in (0,\delta)\\
\var(x), & t=0,\endcases
$$
where $\delta$ is a sufficiently small positive number.

Let $S$ be the Seeley extension operator (see \cite{8}) from
$C^\infty\(\ov{\Bbb
R}^n_+\)$ to $C^\infty\(\Bbb R^{n+1}\)$.

Let $\rho_\delta\in C^\infty_0(-\delta,\delta)$ such that
$\rho_\delta(0)=1$.

Because of Lemma 1.1, Theorem 1.2 will follow from the
following:

\proclaim{Proposition 2.1} Let $\beta=\{Z_1,\cdots,Z_k\}$. If
$V$ and $V'$ are
small enough, $V\subset\subset V'$, there exist $U$ a
neighborhood of $0$ in
$\Bbb R^{n+1}$, $C>0$ and $\delta >0$ such that for any
$\psi\in C^\infty_0(V)$,

$$
\left\|\rho_\delta S(E\psi)\right\|_{W_{1,p}(U)}\le
C\|\psi\|_{W_{1-{1\over p},p}
(V,V',\beta)}.
$$
(Here $\beta=\{Z_1,\cdots,Z_n\}$ is the special basis chosen
in this section).
\endproclaim
\bigskip

The proof of this proposition will be based on some lemmas.
For $\tau=(\tau_1,\cdots,
\tau_n)$ and $x\in V$, let

$$
\eta(\tau,x)=e^{\tau_1Z_1}\cirdot e^{\tau_nZ_n}x
$$
where $\tau=(\tau_1,\cdots,\tau_n)$.

If the neighborhood $V$ of $0$ in $\Bbb R^n_x$ is
sufficiently small, $\eta(\tau,
x)$ is a diffeomorphism from a neighborhood of $0$ in $\tau$
space into $V$.

\proclaim{Lemma 2.2} In the coordinates of $\eta(\tau,x)$, we
have

$$
Z_j={\part\over\part\tau_j}+\sum^n_{\ell=k+1}\vars_{\ell
j}(\tau,x)
{\part\over\part\tau_\ell}
$$
where each $\vars_{\ell j}(\tau,x)=o(\tau)$ for each $1\le
j\le k$;

$$
Z_i=\sum^n_{\ell=k+1}\vars_{i\ell}(\tau,x){\part\over\part\tau_\ell}
\quad\hbox{ for } k+1\le i\le n.
$$
\endproclaim

\demo{Proof} Recall that

$$
Z_i={\part\over\part x_i}+\sum^n_{j=k+1}
a_{ij}(x){\part\over\part x_j}
$$
$\(a_{ij}(x)=a_{ij}(x,0)\)$ for $1\le i\le k$. Moreover, for
each
$\ell\ge k+1$, $\exists\, i,j$ in $\{1,\cdots,k\}$ such that
$Z_\ell=\[Z_i,Z_j\]$.

Therefore, for $\ell\ge k+1$, each

$$
Z_\ell=\sum^n_{j=k+1} b_{\ell j}(x){\part\over\part x_j}
$$
for some smooth  $b_{\ell j}$.

It follows that

$$
\eta(\tau,x)=\(x_1+\tau_1,\cdots,x_k+\tau_k,
x_{k+1}+\tau\cdot
g^{k+1}(\tau,x),\cdots,x
_n+\tau\cdot g^n(\tau,x)\)
$$
where for $i\ge k+1$, $\tau\cdot
g^i(\tau,x)=\sum^n_{\ell=1}\tau_\ell\, g^i_\ell
(\tau,x)$ for some $C^\infty$ functions $g^i_\ell(\tau,x)$.

Moreover, the diffeomorphism $\tau\longmapsto\eta(\tau,x)$
maps each $\disp{
\part\over\part\tau_j}\bigg|_0$ to $Z_j\big|_x$ for $1\le
j\le k$. The
lemma follows from these  remarks.

\enddemo

\proclaim{Lemma 2.3} Let $B_1,\cdots,B_k$ be any operators.
Then

(a) $B_1B_2\cdots B_k-I=(B_1-I)B_2\cdots B_k+(B_2-I)B_3\cdots
B_k+\cdots+(B_k-I)$;

(b) The product $B_1\cdots B_{k-1}$ $(B_k-I)$ is a linear
combination of terms
of the form

$$
\(B_{i_1}-I\)\(B_{i_2}-I\)\cdots \(B_{i_j}-I\),\qquad 1\le
j\le k
$$
and $1\le i_1,<\cdots <i_j\le k$.
\endproclaim
\medskip

Note that (b) can easily be proved by induction and (a) is
obvious. We remark
that (a) was used in \cite{4}.

\demo{Proof of Proposition 2.1} Recall that for $\psi\in
C^\infty_0(V)$ and
$0<t\le\delta$,

$$
H\psi(x,t)=H_1\cirdot H_n\psi(x,t).
$$

We will estimate $\|\part_t
(H\psi)\|_{L^p(V\times(-\delta,\delta))}$ and

$$
\|X_i(H\psi)\|_{L^p(V\times(-\delta,\delta))}.
$$

Observe that $\disp{\part\over\part t}(H\psi)$ is a sum of
terms of the
form

$$
H_1\cirdot \part_t H_i\cirdot H_n\psi,
$$
where for $1\le i\le k$,

$$
\part_t H_i f(x,t)={-1\over t^2}\int^t_0\(e^{\tau Z_i}-I\)
f(x) d\tau +
{1\over t} \(e^{tZ_i}-I\) f(x) \tag 2.1
$$
while for $k+1\le i\le n$,

$$
\part_t H_i f(x,t)={-2\over t^3}\int^{t^2}_0\(e^{\tau
Z_i}-I\) f(x) d\tau +
{2\over t} \(e^{t^2Z_i}-I\) f(x) \tag 2.2
$$

Writing each $H_j$ as $(H_j-I)+I$, we can express
$H_1\cirdot\part_tH_i\cirdot
H_n\psi$ as a sum of  terms of the form:

$$
\(H_{i_1}-I\)\cirdot \part_tH_i\cirdot
\(H_{i_m}-I\)\psi(x,t).
$$

Now if $i_1\le k$, the latter can be bounded by a sum of
terms of the form

$$
{1\over t}\cdot
\left|\psi\(e^{\tau_{i_1}Z_{i_1}}y\)-\psi(y)\right|\tag 2.3
$$
where $y=e^{s_{j_1}Z_{j_1}}\cdots e^{s_{j_\ell}Z_{j_\ell}}x$
for some
$s_{j_1},\cdots, s_{j_\ell}\in [0,t]$.

If $i_1>k$, we use

$$
{1\over t}\cdot
\left|\psi\(e^{\tau^2_{i_1}Z_{i_1}}y\)-\psi(y)\right|,\qquad
\left|\tau_{i_1}\right|\le t. \tag 2.4
$$

Let $V_1\subset\subset V$ be a neighborhood of $0$ such that

$$
e^{\tau_{i_1}Z_{i_1}}\cdots
e^{\tau_{i_m}Z_{i_m}}(V_1)\subseteq V\qquad\hbox{
for }\,\,\left|\tau_{i_j}\right|\le\delta,
$$
and $1\le i_1 <\cdots < i_m\le n$.

From (2.1) $-$ (2.4), we get

$$
\align
& \left\|{\part(H\psi)\over\part
t}\right\|_{L^p\(V_1\times(0,\delta)\)}\\
\le\,\, & C\left\{\sum^k_{i=1}\[\int^\delta_0\(t^{{1\over
p}-1}\omega_i (t,\psi,V)\)
^p{dt\over t}\]^{1\over p}
+\sum^n_{j=k+1}\[\int^\delta_0\(t^{{1\over p}-1}\omega_
j(t^2,\psi,V)\)^p{dt\over t}\]^{1\over p}\right\}.\quad\tag
2.5
\endalign
$$

Next for $j\ge k+1$, we estimate $\omega_j(t^2,\psi,V)$ which
by definition

$$
=\sup_{|\tau|\le
t^2}\left\|\psi\(e^{\tau^2Z_j}x\)-\psi(x)\right\|_{L^p(V)}.
$$

Since $j\ge k+1$, $\,\,\exists\,\, m,\ell\le k$ such that
$Z_j=\[Z_m,Z_\ell\]$.

We have:

$$
\align
\left|\psi\(e^{\tau^2Z_j}x\)-\psi(x)\right|
\le & \left|\psi\(e^{\tau^2Z_j}x\)-\psi\(e^{-\tau
Z_m}e^{-\tau Z_\ell}e^{\tau
Z_m}e^{\tau Z_\ell}x\)\right|\\
&\, + \left|\psi\(e^{-\tau Z_m}e^{-\tau Z_\ell}e^{\tau
Z_m}e^{\tau Z_\ell}x\)
-\psi\(e^{-\tau Z_\ell} e^{\tau Z_m} e^{\tau
Z_\ell}x\)\right|\\
&\, + \left|\psi\(e^{-\tau Z_\ell} e^{\tau Z_m} e^{\tau
Z_\ell} x\)-\psi\(
e^{\tau Z_m} e^{\tau Z_\ell} x\)\right|\\
&\, + \left|\psi\(e^{\tau Z_m} e^{\tau
Z_\ell}x\)-\psi\(e^{\tau Z_\ell} x\)\right|\\
&\, + \left|\psi\(e^{\tau Z_\ell}x\)-\psi (x)\right|.\tag 2.6
\endalign
$$

In the sum on the right in (2.6), every term except the first
one can be estimated
by

$$
\left|\psi\(e^{-\tau Z_m}y\)-\psi(y)\right| +
\left|\psi\(e^{-\tau Z_\ell}y\)-
\psi(y)\right|
$$
where $y$ varies in $V$ provided $x\in V_1$.

To  estimate the $L^p$ norm of the first term, we recall
first from the
Campbell-Hausdorff formula that

$$
\left| e^{\tau^2 Z_j}x - e^{-\tau Z_m} e^{-\tau Z_\ell}
e^{\tau Z_m} e^{\tau Z_
\ell} x\right| = o(\tau^3),
$$
as long as $x$ varies in the relatively compact set $V_1$.

The latter allows us to apply Lemma 3.4 of \cite{4} to
conclude that

$$
\align
\(\int_{V_1} \left|\psi\(e^{\tau^2 Z_j}x\)-\psi \(e^{-\tau
Z_m} e^{-\tau Z_\ell}
e^{\tau Z_m} e^{\tau Z_\ell} x\)\right|^p dx\)^{1\over p}
&\le C\sup_{|s|\le t^3}\left\|\psi (\cdot +s)-\psi
(\cdot)\right\|_{L^p(V)}\\
&= C\,\omega (t^3,\psi,V),
\endalign
$$
where $\omega$ is the usual modulus of $L^p$ continuity.

This inequality together with (2.6) imply that when $j\ge
k+1$,
$$
\omega_j(t^2,\psi,V)\le C\(\sum^k_{i=1}\omega_i(t,\psi,V)
+\omega(t^3,\psi,V)\)
\tag 2.7
$$

Next observe that

$$
\align
\(\int^\delta_0\[t^{{1\over
p}-1}\omega(t^3,\psi,V)\]^p{dt\over t}\)^{1\over p}
&= \(\int^\delta_0\({\omega(s,\psi,V)\over s^{1/3\(1-{1\over
p}\)}}\)^p{ds\over
s}\)^{1\over p}\\
&\le \delta^{1/6\(1-{1\over p}\)}\(\int^\delta_0
\({\omega(s,\psi,V)\over s^{1/2
\(1-{1\over p}\)}}\)^p{ds\over s}\)^{1\over p}\\
&\le C_1 \delta^{1/6\(1-{1\over
p}\)}\left\|H\psi\right\|_{W_{1,p}(V\times
(0,\delta))}
\quad\hbox{ (by (1.5))}.\tag 2.8
\endalign
$$

From (2.5), (2.7) and (2.8)  we conclude:

$$
\align
\|\part_t(H\psi)\|_{L^p(V_1\times (0,\delta))} &\le
C\Bigg\{\sum^k_{i=1}\[\int^\delta
_0\(t^{{1\over p}-1}\omega_i(t,\psi,V)\)^p{dt\over
t}\]^{1\over p}\Bigg.\\
&\,\, +\Bigg.\delta^{1/6\(1-{1\over
p}\)}\|H\psi\|_{W_{1,p}(V\times (0,\delta))}
\Bigg\}.\tag 2.8'
\endalign
$$

We next estimate $\|X_i(H\psi)\|_{L^p(V\times (0,\delta))}$
for $i=1,...,k$.

Recall that $\eta(\tau,x)= e^{\tau_1Z_1}\cdots
e^{\tau_nZ_n}x$ and

$$
H\psi (x,t) = {1\over t^{2n-k}}\int^t_0\cdots\int^{t^2}_0
\psi(\eta(\tau,x))
d\tau\quad (d\tau=d\tau_1\cdots d\tau_n).
$$

We had

$$
\eta_i(\tau,x) =\cases x_i +\tau_i, & i\le k\\
x_i +\tau\cdot g^i(x,t), & i\ge k+1,\endcases
$$
\bigskip

$$
X_j\big|_{t=0} =Z_j = {\part\over\part x_j}+\sum^n_{\ell=k+1}
a_{j\ell}
(x){\part\over\part x_\ell}.
$$

To estimate $X_j(H\psi)$, we will first compute $Z_j(H\psi)$.
We will need to
compare

$$
Z_j\{\psi(\eta(\tau,x))\}\qquad\hbox{ with }\qquad
(Z_j\psi)(\eta(\tau,z)).
$$

From the expressions of the $Z_j$ and the $\eta_i$, we get:

$$
Z_j\{\psi(\eta(\tau,x))\}=(Z_j\psi)(\eta(\tau,x))+\sum^n_{\ell
=k+1} p_{j\ell}
(x,\tau){\part\psi\over\part x_\ell}(\eta(\tau,x)),
$$
where $p_{j\ell}(x,\tau)=o(\tau)$, $ 1\le j\le k$.

Observe that from the form of $\{Z_1,\cdots,Z_n\}$, each
$\disp{\part\over
\part x_\ell}$ for $\ell \ge k+1$ is a linear combination of
$Z_{k+1},
\cdots,Z_n$. Therefore, we get:

$$
Z_j\{\psi(\eta(\tau,x))\}=(Z_j\psi)(\eta(\tau,x))+\sum^n_{\ell
=k+1} q_{j\ell}
(x,\tau)(Z_\ell\psi)(\eta(x,\tau)),\tag 2.9
$$
$q_{j\ell}(x,\tau)=o(\tau)$ and $1\le j\le k$.

By similar arguments, we can get a relation of the type (2.9)
for $Z_i$ when
$i\ge k+1$.

Now, for $1\le i\le k$,

$$
X_i=Z_i+t\sum^n_{j=k+1} C_{ij} (x,t) Z_j\quad\hbox{ for some
} C^\infty C_{ij}.
$$

Hence by (2.9) and its analogue for $j\ge k+1$,

$$
X_i\{\psi(\eta(\tau,x))\}=\(Z_i\psi\)(\eta(\tau,x))+
t\sum^n_{j=k+1} r_{ij}
(x,t,\tau)\(Z_j\psi\)(\eta),
$$
$$
\hbox{for } 1\le i\le k,\quad \hbox{ for some }\quad
C^\infty\,\, r_{ij}.\tag 2.10
$$

Using (2.10), we can write (for $1\le i\le k$)

$$
X_i\{H\psi(x,t)\}= \text{I} + \text{II},
$$
where

$$
\text{I}= {1\over
t^{2n-k}}\int^t_0\cdots\int^{t^2}_0\(Z_i\psi\)(\eta(\tau,x))
d\tau
$$
and

$$
\text{II}=\sum^n_{j=k+1}{1\over
t^{2n-k-1}}\int^t_0\cdots\int^{t^2}_0 r_{ij}(x,t,\tau)
\(Z_j\psi\)(\eta(\tau,x)) d\tau.
$$

Below we will use the notations

$$
\tau_i(a)=(\tau_1,\cdots,\tau_{i-1}, a,
\tau_{i+1},\cdots,\tau_n),
$$

$\Delta_i(s) f(\tau)=f(\tau_i(s))-f(\tau_i(0))$, and the
easily verifiable
identity

$$
\Delta_i(s)\left\{g(\tau)f(\tau)\right\}=\(g(\tau_i(s))\cdot\
Delta_i(s)f(\tau)\)
+\(f(\tau_i(0))\cdot\Delta_i(s)g(\tau)\).
$$

Now by Lemma 2.2, the term

$$
\align
\text{I} = &{1\over
t^{2n-k}}\int^t_0\cdots\int^{t^2}_0\[{\part\over\part
\tau_i}+\sum^n_{\ell=k+1}\vars_{i\ell}(\tau,x){\part\over\part
\tau_\ell}\]
\psi(\eta(\tau,x))d\tau\\
= & {1\over
t^{2n-k}}\int^t_0\cdots\wideh{\int^t_0}\cdots\int^{t^2}_0\Delta_
i(t)
\psi(\eta(\tau,x)) d\tau_1\cdots \wideh{d\tau_i}\cdots
d\tau_n\\
& + {1\over
t^{2n-k}}\sum^n_{\ell=k+1}\int^t_0\cdots\wideh{\int^{t^2}_0}\
cdots
\int^{t^2}_0\Delta_\ell(t^2)\[\vars_{i\ell}(\tau,x)\psi(\eta(
\tau,x))\]d\tau_1
\cdots\wideh{d\tau_\ell}\cdots d\tau_n\\
& - {1\over
t^{2n-k}}\sum^n_{\ell=k+1}\int^t_0\cdots\int^{t^2}_0{\part\vars_
{i\ell}\over\part\tau_\ell}(\tau,x)\psi(\eta(\tau,x))d\tau
\endalign
$$

$$
\align
\quad\,\, = & {1\over
t^{2n-k}}\int^t_0\cdots\wideh{\int^t_0}\cdots\int^{t^2}_0\Delta_i
(t)\psi(\eta(\tau,x)) d\tau_1\cdots\wideh{d\tau_i}\cdots
d\tau_n\\
& + {1\over
t^{2n-k}}\sum^n_{\ell=k+1}\int^t_0\cdots\wideh{\int^{t^2}_0}\
cdots
\int^{t^2}_0
\vars_{i\ell}(\tau_\ell(t^2),x)\Delta_\ell(t^2)\psi(\eta(\tau
,x))
d\tau_1\cdots\wideh{d\tau_\ell}\cdots d\tau_n\\
& + {1\over
t^{2n-k}}\sum^n_{\ell=k+1}\int^t_0\cdots\wideh{\int^{t^2}_0}\
cdots
\int^{t^2}_0\Delta_\ell(t^2)\vars_{i\ell}(\tau,x)\psi(\eta(\tau_
\ell(0),x))
d\tau_1\cdots\wideh{d\tau_\ell}\cdots d\tau_n\\
& - {1\over
t^{2n-k}}\sum^n_{\ell=k+1}\int^t_0\cdots\int^{t^2}_0{\part\vars_
{i\ell}(\tau,x)\over\part\tau_\ell}\psi(\eta(\tau,x))d\tau.
\endalign
$$

Again by Lemma 2.2, a typical term in II is, for some $j\ge
k+1$:

$$
\align
\quad\,\, =& {1\over t^{2n-k-1}}\int^t_0\cdots\int^{t^2}_0
r_{ij}(x,t,\tau)(Z_j\psi)
(\eta(\tau,x)) d\tau\\
 & = {1\over
t^{2n-k-1}}\sum^n_{\ell=k+1}\int^t_0\cdots\int^{t^2}_0 r_{ij}
\vars_{i\ell}(\tau,x){\part\over\part\tau_\ell}\{\psi(\eta(\tau,
x))\}d\tau\\
 & = {1\over
t^{2n-k-1}}\sum^n_{\ell=k+1}\int^t_0\cdots\wideh{\int^{t^
2}_0}
\cdots \int^{t^2}_0
r_{ij}(x,t,\tau_\ell(\tau^2))\vars_{i\ell}(\tau_\ell(t^2),x)
\Delta_\ell(t^2)\psi(\eta)
d\tau_1\cdots\wideh{d\tau_\ell}\cdots d\tau_n\\
 & + {1\over
t^{2n-k-1}}\sum^n_{\ell=k+1}\int^t_0\cdots\wideh{\int^{
t^2}_0}
\cdots \int^{t^2}_0
\Delta_\ell(t^2)\left\{r_{ij}(x,t,\tau)\vars_{i\ell}
(\tau,x)\right\}\psi
(\eta(\tau_\ell(0),x))d\tau_1\cdots\wideh{d\tau_\ell}\cdots
d\tau_n\\
 & - {1\over
t^{2n-k-1}}\sum^n_{\ell=k+1}\int^t_0 \cdots \int^{t^2}_0
\left\{
{\part\over\part\tau_\ell}\(r_{ij}\cdot\vars_{i\ell}\)\right\}
\psi(\eta
(\tau,x))d\tau.
\endalign
$$

By the mean value theorem,

$$
\Delta_\ell(t^2)\vars_{i\ell}(\tau,x)=o(t^2)
=\Delta_\ell(t^2)\{ r_{ij}(x,t,
\tau)\vars_{i\ell}(\tau,x)\},
$$
and recall that for $i\le k$,
$\vars_{i\ell}(\tau,x)=o(\tau)$.

Let  $P_t=\{(\tau_1,\cdots,\tau_n):|\tau_1|\le
t,\cdots,|\tau_k|\le t$ and
$|\tau_i|\le t^2$ for $i>k\}$.

Using the estimates on I and II and the Minkowski inequality
we get:

$$
\align
& \left\|X_i(H\psi(x,t))\right\|_{L^p(V\times(0,\delta))}\\
\le\quad &
C\Bigg[\sum^k_{i=1}\(\int^\delta_0\(t^{-1}\sup_{\tau\in
P_t}\|\Delta_i
(t)\psi(\eta(\tau)\cdot)\|_{L^p(V)}\)^pdt\)^{1\over
p}\Bigg.\\
& +\sum^n_{j=k+1}\(\int^\delta_0\(t^{-1}\sup_{\tau\in
P_t}\|\Delta_j(t^2)\psi
(\eta(\tau)\cdot)\|_{L^p(V)}\)^pdt\)^{1\over p}\\
& +\Bigg.\|\psi\|_{L^p(V)}\Bigg]\\
=\quad & C\(J_1(V)+J_2(V) +\|\psi\|_{L^p(V)}\).\tag 2.11
\endalign
$$

To estimate $J_1(V)$ observe that by using Lemma (2.3) (b) we
have:

$$
\align
& \sup_{\tau\in
P_t}\left\|\Delta_i(t)\psi(\eta(\tau)\cdot)\right\|_{L^p(V)}\\
=\quad & \sup_{\tau\in P_t}\left\|e^{\tau_1Z_1}\cdots
e^{\tau_{i-1}Z_{i-1}}\(e^
{\tau_iZ_i}-I\)e^{\tau_{i+1}Z_{i+1}}\cdots
e^{\tau_nZ_n}\psi(\cdot)\right\|_
{L^p(V)}\\
\le\quad & C\sup_{\tau\in P_t}\left\|e^{\tau_1Z_1}\cdots
e^{\tau_{i-1}Z_{i-1}}
\(e^{\tau_iZ_i}-I\)\psi(\cdot)\right\|_{L^p(V')}\\
\le\quad & C\sum^k_{i=1}\omega_i(t,\psi,V').\tag 2.12
\endalign
$$
where $V'$ is a small neighborhood of $V$. Indeed, by
decreasing $\delta$,
we can make $V'$ as close to $V$ as we wish. It follows that

$$
J_1(V)\le C\sum^k_{i=1}\[\int^\delta_0\(t^{{1\over
p}-1}\omega_i(t,\psi,V')\)^p
{dt\over t}\]^{1\over p}.\tag 2.13
$$

The term $J_2(V)$ can be estimated by using the arguments
employed to establish
(2.7) and (2.8).

This yields:

$$
\align
J_2(V) \le &\;
C\Bigg\{\sum^k_{i=1}\[\int^\delta_0\(t^{{1\over
p}-1}\omega_i(t,
\psi,V')\)^p{dt\over t}\]^{1\over p}\Bigg.\\
& +\Bigg.\delta^{1/6(1-{1\over
p})}\|H\psi\|_{W_{1,p}(V'\times (0,\delta))}\Bigg\}.
\tag 2.14
\endalign
$$

From (2.12), (2.13) and (2.14) we get the following: if
$V_1\subset\subset V$,
then $\delta >0$ can be chosen small enough so that

$$\align
\|X_i(H\psi)\|_{L^p(V\times(0,\delta))} \le &\;
C\Bigg\{\sum^k_{i=1}\[\int^
\delta_0\(t^{{1\over p}-1}\omega_i(t,\psi,V)\)^p{dt\over
t}\]^{1\over p}\Bigg.\\
& +\Bigg.\|\psi\|_{L^p(V)}+\delta^{1/6(1-{1\over
p})}\|H\psi\|_{W_{1,p}(V\times
(0,\delta))}\Bigg\}.\tag 2.14'
\endalign
$$

Observe next that using Minkowski's inequality for integrals
one easily gets:

$$
\|H\psi\|_{L^p(V_1\times(0,\delta))}\le \delta^{1\over
p}\|\psi\|_{L^p(V)}.
\tag 2.15
$$

From (2.8'), (2.14) and (2.15) we get the following: given
$V_1\subset\subset V$
neighborhoods of $0$, there exist $\delta >0$ and $C>0$ such
that

$$
\|H\psi\|_{W_{1,p}(V_1\times (0,\delta))}\le
C\(\|\psi\|_{W_{1-{1\over p},p}
(\beta ,V)}+\delta^{1/6(1-{1\over
p})}\|H\psi\|_{W_{1,p}(V\times (0,\delta))}\).
\tag 2.16
$$

Let now $V_2$ be a neighborhood of $0$ such that

$$
e^{\tau_1Z_1}\cdots e^{\tau_nZ_n}(V_2)\subseteq V_1
$$
for $|\tau_j|\le\delta$.

If $\psi\in C^\infty_0 (V_2)$, then the $x$-support of
$H\psi$ is in $V_1$.

Therefore, for $\psi\in C^\infty_0(V_2)$, the term

$$
\delta^{1/6(1-{1\over
p})}\|H\psi\|_{W_{1,p}(V\times(0,\delta))}
$$
in (2.16) can be absorbed on the left hand side yielding: for
$\psi\in C^\infty
_0(V_2)$,

$$
\|H\psi\|_{W_{1,p}(V_1\times (0,\delta))}\le
C\|\psi\|_{W_{1-{1\over p},p}(\beta,V_1)}.\tag 2.17
$$

The assertion in Proposition 2.1 easily follows from (2.17).

As indicated before, Proposition 2.1 and Lemma 1.1 imply
Theorem 1.2.
\enddemo

\enddocument